\documentclass[a4paper,12pt]{amsart}
\usepackage{amsfonts}
\usepackage{amssymb}
\usepackage{ifthen}
\usepackage{graphicx}
\usepackage{enumitem}
\usepackage{amsmath}
\nonstopmode \numberwithin{equation}{section}
\usepackage[usenames]{color}
\newtheorem{thm}{Theorem}[section]
\newtheorem{lem}{Lemma}[section]
\newtheorem{cor}[thm]{Corollary}
\newtheorem{prop}[thm]{Proposition}

\theoremstyle{definition}
\newtheorem{mlem}{Main lemma}[section]
\newtheorem{assertion}{Assertion}[section]
\newtheorem{cl}{Claim}[section]
\newtheorem{ca}{Case}[section]
\newtheorem{sca}{Subcase}[section]
\newtheorem{scl}{Subclaim}[section]
\newtheorem{conj}[thm]{Conjecture}
\newtheorem{fact}{Fact}[section]
\newtheorem{defn}{Definition}[section]
\newtheorem{prob}{Problem}[section]
\newtheorem{ques}[thm]{Question}
\newtheorem{rem}{Remark}[section]
\newtheorem{exam}{Example}[section]

\numberwithin{equation}{section}
\DeclareMathOperator*{\esssup}{ess\,sup}

\newcounter {own}
\def\theown {\thesection       .\arabic{own}}

\newlist{steps}{enumerate}{1}
\setlist[steps,1]{
  leftmargin=*,
  label=\textbf{Step \arabic*}.,
  ref=Step~\arabic*,
}

\newenvironment{pf}[1][]{%
 \vskip 3mm
 \noindent
 \ifthenelse{\equal{#1}{}}%
  {{\slshape Proof. }}%
  {{\slshape #1.} }%
 }%
{\qed\bigskip}

\newcounter{alphabet}
\newcounter{tmp}

\newcounter{alphabet2}

\makeatletter
\newcommand{\Ref}[1]{\@ifundefined{r@#1}{}{\setcounter{tmp}{\ref{#1}}\Alph{tmp}}}
\makeatother

\newenvironment{Lem}[1][]{\refstepcounter{alphabet}%
\bigskip%
\noindent%
{\bf Lemma \Alph{alphabet}}%
{\bf .} \itshape}{\vskip 8pt}

\newcommand{\ID}{{\mathbb D}}

\newcommand{\real}{{\operatorname{Re}\,}}

\def\be{\begin{equation}}
\def\ee{\end{equation}}

\newcommand{\ben}{\begin{enumerate}}
\newcommand{\een}{\end{enumerate}}

\newcommand{\blem}{\begin{lem}}
\newcommand{\elem}{\end{lem}}
\newcommand{\bthm}{\begin{thm}}
\newcommand{\ethm}{\end{thm}}
\newcommand{\bcor}{\begin{cor}}
\newcommand{\ecor}{\end{cor}}
\newcommand{\beg}{\begin{exam}}
\newcommand{\eeg}{\end{exam}}
\newcommand{\begs}{\begin{examples}}
\newcommand{\eegs}{\end{examples}}
\newcommand{\bdefe}{\begin{defn}}
\newcommand{\edefe}{\end{defn}}
\newcommand{\bprob}{\begin{prob}}
\newcommand{\eprob}{\end{prob}}
\newcommand{\bques}{\begin{ques}}
\newcommand{\eques}{\end{ques}}
\newcommand{\bei}{\begin{itemize}}
\newcommand{\eei}{\end{itemize}}
\newcommand{\bcon}{\begin{conj}}
\newcommand{\econ}{\end{conj}}
\newcommand{\bop}{\begin{op}}
\newcommand{\eop}{\end{op}}

\newcommand{\bas}{\begin{assertion}}
\newcommand{\eas}{\end{assertion}}

\newcommand{\bfa}{\begin{fact}}
\newcommand{\efa}{\end{fact}}

\newcommand{\bca}{\begin{ca}}
\newcommand{\eca}{\end{ca}}
\newcommand{\bsca}{\begin{sca}}
\newcommand{\esca}{\end{sca}}

\newcommand{\bcl}{\begin{cl}}
\newcommand{\ecl}{\end{cl}}

\newcommand{\bmlem}{\begin{mlem}}
\newcommand{\emlem}{\end{mlem}}

\newcommand{\bscl}{\begin{scl}}
\newcommand{\escl}{\end{scl}}

\newcommand{\bcons}{\begin{conjs}}
\newcommand{\econs}{\end{conjs}}

\newcommand{\bprop}{\begin{prop}}
\newcommand{\eprop}{\end{prop}}

\newcommand{\br}{\begin{rem}}
\newcommand{\er}{\end{rem}}
\newcommand{\brs}{\begin{rems}}
\newcommand{\ers}{\end{rems}}
\newcommand{\bo}{\begin{obser}}
\newcommand{\eo}{\end{obser}}
\newcommand{\bos}{\begin{obsers}}
\newcommand{\eos}{\end{obsers}}
\newcommand{\bpf}{\begin{pf}}
\newcommand{\epf}{\end{pf}}
\newcommand{\ba}{\begin{array}}
\newcommand{\ea}{\end{array}}
\newcommand{\beq}{\begin{eqnarray}}
\newcommand{\beqq}{\begin{eqnarray*}}
\newcommand{\eeq}{\end{eqnarray}}
\newcommand{\eeqq}{\end{eqnarray*}}

\newcounter{minutes}
\divide\time by 60
\newcounter{hours}
\multiply\time by 60 \addtocounter{minutes}{-\time}

\begin{document}

\bibliographystyle{amsplain}

\title{Norm estimates of the Cauchy transform and related operators}

\author{Jian-Feng Zhu}
\address{Jian-Feng Zhu, School of Mathematical Sciences,
Huaqiao University,
Quanzhou 362021, People's Republic of China and Department of Mathematics,
Shantou University, Shantou, Guangdong 515063, People's Republic of China.
} \email{flandy@hqu.edu.cn}

\author{David Kalaj}
\address{David Kalaj, University of Montenegro, Faculty of Natural Sciences and
Mathematics, Cetinjski put b.b. 81000 Podgorica, Montenegro}
\email{davidk@ac.me}

\date{}

\begin{abstract}
Suppose $f\in L^p(\mathbb{D})$, where $p\geq1$ and $\mathbb{D}$ is the unit disk.
Let $\mathfrak{J}_0$ be the integral operator defined as follows:
$\mathfrak{J}_0[f](z)=\int_{\mathbb{D}}\frac{z}{1-\bar{w}z}f(w)\mathrm{d}A(w)$, where $z$, $w\in\mathbb{D}$ and $\mathrm{d}A(w)=\frac{1}{\pi}\mathrm{d}x\mathrm{d}y$ is the normalized area measure on $\mathbb{D}$.
Suppose $\mathfrak{J}_0^*$ is the adjoint operator of $\mathfrak{J}_0$. Then $\mathfrak{J}^*_0=\mathfrak{B}\mathfrak{C}$, where $\mathfrak{B}$ and $\mathfrak{C}$ are the operators induced by the Bergman projection and Cauchy transform, respectively.
In this paper, we obtain the $L^1$, $L^2$ and $L^{\infty}$ norm of the operator $\mathfrak{J}_0^*$.
Moreover, we obtain the $L^p(\mathbb{D})\rightarrow L^\infty(\mathbb{D})$ norm of the operators $\mathfrak{C}$ and $\mathfrak{J}_0^*$, provided that $p>2$.
This study is a continuation of the investigations carried out in \cite{Baranov} and \cite{kalaj1}.
\end{abstract}
\thanks{}

\maketitle \pagestyle{myheadings} \markboth{Jian-Feng Zhu and David Kalaj}{Norm estimates of the Cauchy transform and related operators}


\section{Introduction}\label{sec-1}
Let $\ID=\{z:|z|<1\}$ be the unit disk of $\mathbb{C}$. Denote by $L^p(\ID) (1\leq p\leq\infty)$ the space of
complex-valued measurable functions on $\ID$ with finite integral
$$\|f\|_{p}=\left(\int_{\ID}|f(z)|^p\mathrm{d} A(z)\right)^{\frac{1}{p}},\ \ \ 1\leq p<\infty,$$
where
$$\mathrm{d}A(z) =\frac{1}{\pi}\mathrm{d}x\mathrm{d}y=\frac{1}{\pi}r\mathrm{d}r\mathrm{d}\theta,\ \ \ z=x+iy=re^{i\theta},$$
is the normalized area measure on $\ID$ (cf. \cite[Page 1]{Hendenmalm}). For the case $p=\infty$, we let $L^{\infty}(\ID)$ denote the space of (essentially) bounded
functions on $\ID$. For $f\in L^{\infty}(\ID)$, we define
$$\|f\|_\infty=\esssup\{|f(z)|: z\in\ID\}.$$
The space $L^{\infty}(\ID)$ is a Banach space with the above norm  (cf. \cite[Page 2]{Hendenmalm}).
\subsection{The Cauchy transform $\mathfrak{C}$}
Let $\Omega\subset\mathbb{C}$ be a bounded domain in the complex plane.
It follows from \cite[Page 7]{Baranov} that the {\it Cauchy integral operator} ({\it Cauchy transform}) $\mathfrak{C}_{\Omega}: L^p(\Omega)\rightarrow L^{p}(\Omega)$ is defined by (see also \cite{Anderson, kalaj1})
$$\mathfrak{C}_{\Omega}[f](z)=\int_{\Omega}\frac{f(w)}{w-z}\mathrm{d} A(w).$$
Unlike the two-dimensional {\it Hilbert transform} $\mathcal{H}f$ (also called the {\it Beurling transform}), the Cauchy transform is not bounded as an
operator from $L^2(\mathbb{C})$ to $L^2(\mathbb{C})$. The characteristic function of any bounded domain has Cauchy transform whose modulus behaves like
$|\zeta|^{-1}$ as $\zeta\rightarrow\infty$, and so does not belong to $L^2(\mathbb{C})$ (cf. \cite{Anderson}).
Nonetheless, it follows from the Sobolev embedding theorem (cf. \cite{Anderson, Mazja}) that $\mathfrak{C}_{\Omega}$ is a bounded operator from $L^2(\Omega)$ to $L^p(\Omega)$ for all $p<\infty$. More precisely, we observe that if $f\in L^2(\Omega)$, then
$$\frac{\partial}{\partial\bar{z}}\mathfrak{C}_{\Omega}[f](z)=-f(z)\in L^2(\mathbb{C}),$$
and (cf. \cite[Page 7]{Baranov}, see also \cite[Page 157, (7.10)]{Lehto})
$$\frac{\partial}{\partial z}\mathfrak{C}_{\Omega}[f](z)=\mathcal{H}f(z)\in L^2(\mathbb{C}).$$

Throughout this paper, we consider the case $\Omega=\ID$.
For simplicity, we write $\mathfrak{C}$ instead of $\mathfrak{C}_{\ID}$ for the Cauchy transform of $\ID$.

Recall that the norm of an operator $T: X\rightarrow Y$ between normed spaces $X$ and $Y$ is defined by
$$\|T\|_{X\rightarrow Y}=\sup\{\|Tx\|_Y:\|x\|_X=1\}.$$
For the case of $X=Y=L^{p}(\Omega)$, we write $\|T\|_p$ instead of $\|T\|_{L^{p}(\Omega)\rightarrow L^{p}(\Omega)}$ for
the $L^p$ norm of the operator $T$.

The norm estimates of the Cauchy transform on $L^p(\Omega)$, i.e.,
$\|\mathfrak{C}_{\Omega}\|_{p}$,
has been studied by many mathematicians, but is still not known.
In the case of $\Omega=\ID$, it was proved in \cite{Anderson} that
$$\|\mathfrak{C}\|_{2}=\frac{2}{j_0},$$
where $j_0\approx2.4048256$ is the smallest positive zero of the Bessel function $J_0=\sum\limits_{k=0}^{\infty}\frac{(-1)^k}{(k!)^2}\left(\frac{x}{2}\right)^{2k}$.

Since $\|\mathfrak{C}\|_{1}=2$ (cf. \cite[(12)]{Dostanic}), using the {\it Riesz-Thorin interpolation theorem} (cf. \cite[Theorem 1.1.1]{Bergh}), it was proved in \cite[Theorem 1]{Dostanic} that
$$\|\mathfrak{C}\|_{p}\leq 2\cdot j_0^{-2(1-\frac{1}{p})}, \ \ \ \mbox{for}\ \ \ 1\leq p\leq2$$
and
$$\|\mathfrak{C}\|_{p}\leq 2\cdot j_0^{-\frac{2}{p}}, \ \ \ \mbox{for}\ \ \ p\geq 2.$$

\subsection{The related operator $\mathfrak{J}_0$ and its adjoint operator $\mathfrak{J}_0^*$.}
Let $\mathfrak{J}_0: L^p(\ID)\rightarrow L^p(\ID)$ be the integral operator defined by (cf. \cite[Page 9 and Page 12]{Baranov})
$$\mathfrak{J}_0[f](z)=\int_{\ID}\frac{z}{1-\bar{w}z}f(w)\mathrm{d}A(w) .$$
Suppose $\mathfrak{J}_0^*$ is the adjoint operator of $\mathfrak{J}_0$. The following transform
$$\mathfrak{B}[f](z)=\int_{\ID}\frac{f(w)}{(1-\bar{w}z)^2}\mathrm{d}A(w) $$
is the Bergman projection (here we refer to \cite{Hendenmalm,Perala,Zhukh} and the references therein for more discussions of Bergman projection).
It was proved in \cite[Page 12]{Baranov} that $\mathfrak{J}_0^*=\mathfrak{BC}$, the composition of $\mathfrak{B}$ and $\mathfrak{C}$,
where
$$\mathfrak{J}_0^*[f](z)=\int_{\ID}\frac{\bar{w}}{1-\bar{w}z}f(w)\mathrm{d}A(w) .$$

Another relation of the Cauchy transform $\mathfrak{C}$ and the related operator $\mathfrak{J}_0^*$ was given in \cite{Astala} and \cite{kalaj1} as follows:
For $f\in L^p(\ID)$, where $1<p<\infty$, {\it the Cauchy transform for Dirichlet's problem} (see \cite[Page 155]{Astala}) of $f$ is defined by
$$C_{\Delta}[f](z)=\int_{\ID}\left(\frac{1}{z-w}+\frac{\bar{w}}{1-\bar{w}z}\right)f(w)\mathrm{d}A(w) , \ \ \ z\in\ID.$$
The operator $C_{\Delta}$ is hence induced by the $z$-derivative of the {\it Green's function} and $\frac{\partial}{\partial{\bar{z}}}C_{\Delta}[f]=f$.
Obviously, one has
$$C_{\Delta}[f](z)=\mathfrak{J}_0^*[f](z)-\mathfrak{C}[f](z).$$
We refer to \cite{Baranov} for more discussions on the relations of the operators $\mathfrak{C}$, $\mathcal{H}f$ and $\mathfrak{J}_0^*$.

A related result, but for the so-called Cauchy operator for the Dirichlet problem in the unit disk, has been given by the second author of this paper in \cite{kalaj1,kalaj2}.
Precisely, it was proved in \cite[Theorem A]{kalaj1} that:
$$\|C_{\Delta}\|_{p}\leq 2\cdot j_0^{-2(1-\frac{1}{p})}, \ \ \ \mbox{for}\ \ \ 1\leq p\leq2$$
and
$$\|C_{\Delta}\|_{p}\leq \frac{4}{3}\left(\frac{2j_0}{3}\right)^{-\frac{2}{p}}, \ \ \ \mbox{for}\ \ \ p\geq 2.$$
The equalities of the above inequalities can be attained for the case $p=1$, $p=2$ and $p=\infty$.

In this paper, we obtain the following norm estimates: $\|\mathfrak{C}\|_{L^p(\ID)\rightarrow L^{\infty}(\ID)}$,
$\|\mathfrak{J}_0\|_{L^p(\ID)\rightarrow L^{\infty}(\ID)}$ and $\|\mathfrak{J}_0^*\|_{L^p(\ID)\rightarrow L^{\infty}(\ID)}$,  for $p>2$.
Furthermore, we obtain the following norms: $\|\mathfrak{J}_0^*\|_{1}$, $\|\mathfrak{J}_0^*\|_{2}$ and $\|\mathfrak{J}_0^*\|_{\infty}$. Our main results are as follows.

\subsection{The norm estimates of operators from $L^p(\ID)$ to $L^{\infty}(\ID)$.}

\begin{thm}\label{2019-May-26-thm}
For $p>2$,
\be\label{2019-May-20-1}
\|\mathfrak{C}\|_{L^p(\ID)\rightarrow L^{\infty}(\ID)}=\left(\frac{2p-2}{p-2}\right)^{1-\frac{1}{p}}.
\ee
If in particular $p=\infty$, then $$\|\mathfrak{C}\|_{\infty}=2.$$
\end{thm}

\begin{thm}\label{thm-Lp-BC*}
For $p>2$,
\be\label{2019-June-30-1}
\|\mathfrak{J}_0\|_{L^p{(\ID)}\rightarrow L^{\infty}(\ID)}=\left(\frac{\Gamma(\frac{p-2}{p-1})}{\Gamma^2\left(\frac{3p-4}{2p-2}\right)}\right)^{1-\frac{1}{p}},
\ee
where $\Gamma$ is the Gamma function.

If in particular $p=\infty$, then
\be\label{June-29-4}\|\mathfrak{J}_0\|_{\infty}=\frac{4}{\pi}.\ee
\end{thm}

\begin{thm}\label{2019-June-11-thm}
For $p>2$,
\be\label{2019-June-11-1}
\|\mathfrak{J}_0^*\|_{L^p{(\ID)}\rightarrow L^{\infty}(\ID)}= A(p)^{1-\frac{1}{p}},
\ee
where
$$A(p)=2\frac{{}_3F_2[1+\frac{p}{2(p-1)}, \frac{p}{2(p-1)},\frac{p}{2(p-1)}; 1, 2+\frac{p}{2(p-1)}; 1 ]}{2+\frac{p}{p-1}}$$
and ${}_3F_2$ is the hypergeometric function given by $(\ref{hypergeometric})$.

If in particular $p=\infty$, then
\be\label{June-29-1}\|\mathfrak{J}_0^*\|_{\infty}=\frac{1+2\alpha}{\pi},\ee
where $\alpha\approx 0.915966$ is the Catalan's constant.
\end{thm}

\begin{rem}\label{Rem-3.2}
(1)
For $p>2$, let $q=\frac{p}{p-1}\in[1, 2)$. According to the definition of the hypergeometric function, we see that
$$\frac{{}_3F_2[1+\frac{q}{2}, \frac{q}{2},\frac{q}{2}; 1, 2+\frac{q}{2}; 1 ]}{2+q}=\sum\limits_{n=0}^{\infty}\left(\frac{\Gamma(n+\frac{q}{2})}{n!\Gamma(\frac{q}{2})}\right)^2\frac{1}{2n+q+2}.$$
Using Lemma \Ref{DZ-2019-lem-0} below, we have: for $n\geq1$,
$$\frac{\Gamma(n+\frac{q}{2})}{n!}\leq\frac{1}{n^{1-\frac{q}{2}}}.$$
Thus,
\begin{align*}
A(p)&\leq2\left(\frac{1}{2+q}+\frac{1}{\Gamma^2(\frac{q}{2})}\sum\limits_{n=1}^{\infty}\frac{1}{n^{2-q}}\frac{1}{2n+q+2}\right)\\ \nonumber
&<2\left(\frac{1}{2+q}+\frac{\zeta(3-q)}{2\Gamma^2(\frac{q}{2})}\right),
\end{align*}
where $\zeta$ is Riemann's zeta function. This shows that $A(p)$ is finite for any $p>2$.

(2) For the case of $1\leq p\leq2$, we show in Remark \ref{rem1}, Remark \ref{Remark-3.2-July} and Remark \ref{rem3.2} that the operators
$\mathfrak{C}$, $\mathfrak{J}_0$ and $\mathfrak{J}_0^*$ will not send $L^p(\ID)$ to $L^{\infty}(\ID)$, respectively.
\end{rem}

\subsection{The $L^1$ norm and $L^2$ norm of $\mathfrak{J}_0^*$.}
The following Corollary \ref{thm-L1-BC} easily follows from (\ref{June-29-4}), since
the $L^1$ norm of an operator is equal to the $L^\infty$ norm of its adjoint operator.

\begin{cor}\label{thm-L1-BC}
\be\label{June-29-L1}\|\mathfrak{J}_0^*\|_{1}=\frac{4}{\pi}.\ee
\end{cor}

For the $L^2$ norm we have:

\begin{thm}\label{thm-L2-BC}
\be\label{June-29-2}\|\mathfrak{J}_0^*\|_{2}^2=\frac{1}{2}.\ee
\end{thm}
Using the Riesz-Thorin interpolation theorem (\cite[Theorem 1.1.1]{Bergh}) together with (\ref{June-29-1}), (\ref{June-29-L1}) and (\ref{June-29-2}), we have the following corollary.
\begin{cor}\label{Reiz-Thorin}
$$\|\mathfrak{J}_0^*\|_{p}\leq \left(\frac{1}{2}\right)^{\frac{1}{p}}\left(\frac{1+2\alpha}{\pi}\right)^{1-\frac{2}{p}},\ \ \ \mbox{for}\ \ \ p\geq2,$$
and
$$\|\mathfrak{J}_0^*\|_{p}\leq \left(\frac{1}{2}\right)^{1-\frac{1}{p}}\left(\frac{4}{\pi}\right)^{\frac{2}{p}-1},\ \ \ \mbox{for}\ \ \ 1\leq p\leq2.$$
The equalities can be attained in the above inequalities for $p=1$, $p=2$ and $p=\infty$.
\end{cor}

The proofs of the above theorems are given in Section \ref{sec-3}.

\section{Preliminaries}\label{sec-2}
In this section, we should recall some known results, and prove three useful lemmas, and one proposition.

\begin{defn} $($cf. \cite[(2.1.2)]{Landrews}$)$
The hypergeometric function
$${}_pF_q[a_1, a_2, \ldots, a_p; b_1, b_2, \ldots, b_q; x]$$
is defined by the series
\be\label{hypergeometric}{}_pF_q[a_1, a_2, \ldots, a_p; b_1, b_2, \ldots, b_q; x]=\sum\limits_{n=0}^{\infty}\frac{(a_1)_n\cdots(a_p)_n}{(b_1)_n\cdots(b_q)_n}\frac{x^n}{n!}\ee
for all $|x|<1$ and by continuation elsewhere.
\end{defn}

Here $(q)_n$ is the Pochhammer symbol which is defined as follows
$$(q)_n=\left\{
\begin{array}
{r@{\ }l}
1, \ \  & \mbox{if}\ \ \ n=0;\\
\\
q(q+1)\cdots(q+n-1), \ \    & \mbox{if}\ \ \  n>0.
\end{array}\right.$$

\begin{Lem}$($\emph{cf.} \cite[Theorem 2.2.2]{Landrews}$)$\label{DZ-2019-lemma-B}
$${}_2F_1(a, b; c; 1)=\frac{\Gamma(c)\Gamma(c-a-b)}{\Gamma(c-a)\Gamma(c-b)}, \ \ \ \mbox{if}\ \ \ \real{(c-a-b)}>0,$$
where $\Gamma$ is the Gamma function.
\end{Lem}

Let $q\in[1, 2)$. The following equality easily follows from Lemma \Ref{DZ-2019-lemma-B}
$${}_2F_1\left[\frac{q}{2}, \frac{q}{2}; 2; 1\right]=\frac{\Gamma(2-q)}{\Gamma^2(2-\frac{q}{2})}.$$

\begin{Lem}$($\emph{cf.} \cite[(7)]{Gautschi}$)$\label{DZ-2019-lem-0}
For $1\leq q\leq2$ and $n=1, 2, \cdots$,
$$\frac{1}{(n+1)^{1-\frac{q}{2}}}\leq\frac{\Gamma(n+\frac{q}{2})}{n!}\leq\frac{1}{n^{1-\frac{q}{2}}},$$
where $\Gamma$ is the Gamma function.
\end{Lem}

\begin{lem}\label{DZ-2019-lem-1}
For $1\leq q<2$, let
$$F(t)=\left(1-t\right)^{2-q} {}_2F_1\left[1-\frac{q}{2}, 2-\frac{q}{2}; 1; t\right],$$
where ${}_2F_1$ is the hypergeometric function.
Then $F(t)$ is a decreasing function of $t$ for $0\leq t\leq1$.
\end{lem}
\bpf
Elementary calculations lead to (see \cite{Landrews})
$$\frac{d}{dt}\left({}_2F_1\left[1-\frac{q}{2}, 2-\frac{q}{2}; 1; t\right]\right)=\frac{1}{2}(2-q)\left(2-\frac{q}{2}\right){}_2F_1\left[2-\frac{q}{2}, 3-\frac{q}{2}; 2; t\right].$$
Then
\be\label{F-derivative}F'(t)=\frac{1}{2}(q-2)(1-t)^{1-q}H(t),\ee
where
$$H(t)=2{}_2F_1\left[1-\frac{q}{2}, 2-\frac{q}{2}; 1; t\right]-\left(2-\frac{q}{2}\right)(1-t){}_2F_1\left[2-\frac{q}{2}, 3-\frac{q}{2}; 2; t\right].$$
Using the definition of the hypergeometric function, we have the following power series:
$$H(t)=\sum\limits_{m=0}^{\infty}a_mt^m,$$
where the coefficients are as follows:
$$a_m=-\frac{\Gamma(1+m-\frac{q}{2})\Gamma(2+m-\frac{q}{2})}{\Gamma(1+m)\Gamma(2+m)\Gamma(2-\frac{q}{2})\Gamma(-\frac{q}{2})}, \ \ \ (m=0, 1, 2, \cdots).$$
For any $1\leq q<2$, we have
$$\Gamma(-\frac{q}{2})<0,$$
and thus, $a_m\geq0$ for all $m\geq0$. This shows that $H(t)\geq0$.
It follows from $(\ref{F-derivative})$ and the assumption $1\leq q<2$ that
$$F'(t)\leq0.$$
Therefore, we know that $F(t)$ is a decreasing function of $t$.

The proof of Lemma $\ref{DZ-2019-lem-1}$ is complete.
\epf

The following result is useful and will be used in proving our main theorems.
For $\beta>0$, $z\in\ID$ and $\zeta=e^{i\theta}\in\mathbb{T}$, where $\mathbb{T}$ is  the unit circle, we have
$$\frac{1}{(1-z\zeta)^\beta}=\sum\limits_{n=0}^{\infty}\frac{\Gamma(n+\beta)}{n!\Gamma(\beta)}z^n\zeta^n.$$
Using Parseval's theorem, one gets
\be\label{2019-May-26-P}\frac{1}{2\pi}\int_0^{2\pi}\frac{d\theta}{|1-ze^{i\theta}|^{2\beta}}=\sum\limits_{n=0}^{\infty}\left(\frac{\Gamma(n+\beta)}{n!\Gamma(\beta)}\right)^2|z|^{2n}.\ee

\begin{lem}\label{DZ-2019-lem-I1(z)}
For $z\in\ID$, let
$$I_1(|z|)=\int_{\ID}\frac{|w|}{|1-\bar{w}z|}\mathrm{d}A(w) .$$
Then
$$\sup\limits_{z\in\ID}I_1(|z|)= \frac{1+2\alpha}{\pi},$$
where $\alpha\approx 0.915966$ is the Catalan's constant.
\end{lem}
\bpf
Let $w=re^{it}\in\ID$. It follows from (\ref{2019-May-26-P}) that
\begin{align*}
I_1(|z|)&=\frac{1}{\pi}\int_0^1r^2\mathrm{d} r\int_0^{2\pi}\frac{1}{|1-\bar{z}re^{it}|}\mathrm{d} t\\ \nonumber
&=2\sum\limits_{n=0}^{\infty}\left(\frac{\Gamma(n+\frac{1}{2})}{n!\Gamma(\frac{1}{2})}\right)^2\frac{1}{2n+3}|z|^{2n}.
\end{align*}
This implies that $I_1(|z|)$ is an increasing function of $|z|$, and has its supremum $I_1(1)$.
Using the equality
$$\sum\limits_{n=0}^{\infty}\left(\frac{\Gamma(n+\frac{1}{2})}{n!\Gamma(\frac{1}{2})}\right)^2\frac{1}{2n+3}=\frac{1+2\alpha}{2\pi},$$
 we have
$$\sup\limits_{z\in\ID}I_1(|z|)= \frac{1+2\alpha}{\pi}.$$

This completes the proof of Lemma \ref{DZ-2019-lem-I1(z)}.
\epf

\begin{lem}\label{DZ-2019-lem-I1(w)}
For $w\in\ID$, let $$I_2(|w|)=\int_{\ID}\frac{|w|}{|1-\bar{w}z|}\mathrm{d}A(z) .$$
Then
$$\sup\limits_{w\in\ID}I_2(|w|)= \frac{4}{\pi}.$$
\end{lem}
\bpf
Following the proof of Lemma \ref{DZ-2019-lem-I1(z)}, we have for $z=re^{it}\in\ID$,
\begin{align*}
I_2(|w|)&=\frac{|w|}{\pi}\int_0^1r\mathrm{d} r\int_0^{2\pi}\frac{1}{|1-\bar{w}re^{it}|}\mathrm{d} t\\ \nonumber
&=|w|\sum\limits_{n=0}^{\infty}\left(\frac{\Gamma(n+\frac{1}{2})}{n!\Gamma(\frac{1}{2})}\right)^2\frac{1}{n+1}|w|^{2n}\\ \nonumber
&={}_2F_1\bigg[\frac{1}{2}, \frac{1}{2}; 2; |w|^2\bigg]|w|.
\end{align*}
This implies that $I_2(|w|)$ is an increasing function of $|w|$. Using the equality
\be\label{2dpi}{}_2F_1\bigg[\frac{1}{2}, \frac{1}{2}; 2; 1\bigg]=\frac{4}{\pi},\ee
we have
$$\sup\limits_{w\in\ID}I_2(|w|)= \frac{4}{\pi}.$$

The proof of Lemma \ref{DZ-2019-lem-I1(w)} is complete.
\epf

\begin{prop}\label{thm-Lp-BC}
$\|\mathfrak{J}_0^*\|_p$ is finite for any $p\geq1$.
\end{prop}
\bpf
According to Lemma \ref{DZ-2019-lem-I1(z)}, we see that for every $z\in\ID$,
$$d\mu(w)=\frac{|w|}{|1-z\bar{w}|}\frac{\mathrm{d}A(w) }{I_1(|z|)}$$
is a probability measure in the unit disk, i.e., $\int_{\ID}d\mu(w)=1$.
For any $p\geq1$ and $f\in L^p(\ID)$, one has
$$|\mathfrak{J}_0^*[f](z)|^p\leq\left(I_1(|z|)\right)^p\left(\int_{\ID}\frac{|w|}{|1-z\bar{w}|}|f(w)|\frac{\mathrm{d}A(w) }{I_1(|z|)}\right)^p.$$
Using Jensen's inequality, we have
$$|\mathfrak{J}_0^*[f](z)|^p\leq I_1^{p-1}(|z|)\int_{\ID}\frac{|w|}{|1-z\bar{w}|}|f(w)|^p\mathrm{d}A(w) .$$
Observe that, since $f\in L^p(\ID)$, it follows from Lemma \ref{DZ-2019-lem-I1(z)} and Lemma \ref{DZ-2019-lem-I1(w)} that
$$\frac{|w|}{|1-z\bar{w}|}|f(w)|^p\in L^1(\ID\times \ID).$$
Using Fubini's theorem, we obtain that
\begin{align*}
\int_{\ID}|\mathfrak{J}_0^*[f](z)|^p\mathrm{d}A(z) &\leq\int_{\ID}I_1^{p-1}(|z|)\mathrm{d}A(z) \int_{\ID}\frac{|w|}{|1-z\bar{w}|}|f(w)|^p\mathrm{d}A(w) \\ \nonumber
&=\int_{\ID}|f(w)|^p\mathrm{d}A(w) \int_{\ID}I_1^{p-1}(|z|)\frac{|w|}{|1-z\bar{w}|}\mathrm{d}A(z) .
\end{align*}

On the other hand, it follows from Lemma \ref{DZ-2019-lem-I1(z)} that
$$I_1(|z|)\leq \frac{1+2\alpha}{\pi},$$
where $\alpha\approx 0.915966$ is the Catalan's constant.
Also, Lemma \ref{DZ-2019-lem-I1(w)} shows that
$$\int_{\ID}\frac{|w|}{|1-z\bar{w}|}\mathrm{d}A(z) \leq \frac{4}{\pi}.$$
Then,
$$\int_{\ID}|\mathfrak{J}_0^*[f](z)|^p\mathrm{d}A(z) \leq\frac{4}{\pi^p}(1+2\alpha)^{p-1}\int_{\ID}|f(w)|^p\mathrm{d}A(w) ,$$
which shows that
$$\|\mathfrak{J}_0^*[f]\|_{p}\leq \frac{4^{\frac{1}{p}}}{\pi}(1+2\alpha)^{1-\frac{1}{p}}\|f\|_{p}.$$

The proof of Proposition \ref{thm-Lp-BC} is complete.
\epf
\section{Proofs of the main results}\label{sec-3}
\subsection*{Proof of Theorem \ref{2019-May-26-thm}}
For $f\in L^p(\ID)$ and $q\in\mathbb{R}$ such that $\frac{1}{p}+\frac{1}{q}=1$, using H\"older's inequality for integrals, we see that
\begin{align}\label{May-26-thm-1}
|\mathfrak{C}[f](z)| &\leq \left(\int_{\ID}|f(w)|^p\mathrm{d}A(w) \right)^{1/p}\left(\int_{\ID}\frac{1}{|w-z|^q}\mathrm{d}A(w) \right)^{1/q}\\ \nonumber
 &=\|f\|_{p}\bigg(\int_{\ID}\frac{1}{|w-z|^{\frac{p}{p-1}}}\mathrm{d}A(w) \bigg)^{1-1/p}.
\end{align}
The assumption of $p>2$ ensures that $q=\frac{p}{p-1}\in[1, 2)$. Let
\be\label{Iq}
K_p(|z|)=\int_{\ID}\frac{1}{|w-z|^{\frac{p}{p-1}}}\mathrm{d}A(w) , \ \ \ \ \ \ z\in\ID.
\ee
We first estimate $K_p(|z|)$ as follows: Using the M\"obius transformation $w=\frac{z-a}{1-\bar{z}a}$, where $a=re^{i\theta}\in \ID$, and
the following equality which comes from (\ref{2019-May-26-P}):
$$\frac{1}{2\pi}\int_0^{2\pi}\frac{\mathrm{d}\theta}{|1-\bar{z}re^{i\theta}|^{4-q}}=\sum\limits_{n=0}^{\infty}\left(\frac{\Gamma(n+2-q/2)}{n!\Gamma(2-q/2)}\right)^2|rz|^{2n},$$
we have
\begin{align*}
K_p(|z|) &= \frac{(1-|z|^2)^{2-q}}{\pi}\int_0^1r^{1-q}\mathrm{d}r\int_0^{2\pi}\frac{\mathrm{d}\theta}{|1-\bar{z}re^{i\theta}|^{4-q}}\\
 &=2 (1-|z|^2)^{2-q}\sum\limits_{n=0}^{\infty}\left(\frac{\Gamma(n+2-q/2)}{n!\Gamma(2-q/2)}\right)^2\frac{|z|^{2n}}{2n+2-q}\\
 &= \frac{2  \left(1-|z|^2\right)^{2-q} {}_2F_1\left[1-\frac{q}{2}, 2-\frac{q}{2}; 1; |z|^2\right]}{2-q},
\end{align*}
where ${}_2F_1$ is the hypergeometric function.
Since $1\leq q<2$, it follows from Lemma $\ref{DZ-2019-lem-1}$ that
$$F(t)=\left(1-t\right)^{2-q} {}_2F_1\left[1-\frac{q}{2}, 2-\frac{q}{2}; 1; t\right]$$
is a decreasing function of $t$ for $0\leq t\leq1$. Therefore, we have  $K_p(|z|)$ is a decreasing function of $|z|$ and
has its maximum
$$K_p(0)= \frac{2}{2-q}.$$
Then using (\ref{May-26-thm-1}), we have
$$\|\mathfrak{C}[f]\|_{\infty}\leq\|f\|_p K_p^{1-\frac{1}{p}}(0),$$
which implies that
\be\label{Kq}\|\mathfrak{C}\|_{\infty}\leq K_p^{1-\frac{1}{p}}(0).\ee

To show the equality of (\ref{Kq}), fix $b\in\ID$ and consider the function
$$f(w)=K_p^{-\frac{1}{p}}(|b|)\frac{w-b}{|w-b|^{\frac{p}{p-1}}},$$
where $p>2$ and $q=\frac{p}{p-1}\in[1, 2)$.
Then
$$\int_{\ID}|f(w)|^p\mathrm{d}A(w)=K_p^{-1}(|b|)\int_{\ID}\frac{1}{|w-b|^q}\mathrm{d}A(w)=1.$$
This shows that $\|f\|_p^p=1$.
On the other hand, elementary calculations show that
$$\left|\mathfrak{C}[f](b)\right|=K_p^{-\frac{1}{p}}(|b|)\int_{\ID}\frac{1}{|w-b|^q}\mathrm{d}A(w)=K_p^{1-\frac{1}{p}}(|b|).$$
Hence,
$$\|\mathfrak{C}\|_{\infty}\geq\|\mathfrak{C}[f]\|_{\infty}\geq K_p^{1-\frac{1}{p}}(b).$$
Observe that since $b$ can be arbitrarily close to $0$, we have
\be\label{2020-Jan-10}\|\mathfrak{C}\|_{\infty}\geq K_p^{1-\frac{1}{p}}(0).\ee
According to (\ref{Kq}) and (\ref{2020-Jan-10}), we see that (\ref{2019-May-20-1}) holds true.

The proof of Theorem \ref{2019-May-26-thm} is complete. \qed
\begin{rem}\label{rem1}
If $1\leq p\leq2$, then $\mathfrak{C}$ will not send $L^p(\ID)$ to $L^{\infty}(\ID)$.
We only need to consider the case of $p=2$, because $L^2(\ID)\subseteq L^p(\ID)$.

Fix $b\in\ID$ and consider the function
$$g(w)=\frac{1}{(\bar{b}-\bar{w})\log\frac{3}{|b-w|}}.$$
Then
$$\int_{\ID}|g(w)|^2\mathrm{d}A(w) =\int_{\ID}\frac{\mathrm{d}A(w) }{|b-w|^2\log^2\frac{3}{|b-w|}}.$$
Let $\xi=b-w$ and $\ID'=\{\xi: |\xi-b|<1\}$. We have $\ID'\subset\ID(0, 2):=\{\xi: |\xi|<2\}$. Thus, for $\xi=R e^{i\theta}\in\ID'$, we have
$$\int_{\ID}|g(w)|^2\mathrm{d}A(w) =\int_{\ID'}\frac{\mathrm{d}A(\xi)}{|\xi|^2\log^2\frac{3}{|\xi|}}\leq\frac{1}{\pi}\int_0^{2\pi}\mathrm{d}\theta\int_0^2\frac{\mathrm{d}R}{R\log^2\frac{3}{R}}=\frac{2}{\log\frac{3}{2}},$$
which shows $g\in L^2(\ID)$.

However, let $D(b)=\{w:|w-b|<1-|b|\}\subset\ID$. Then
$$|\mathfrak{C}[g](b)|=\int_{\ID}\frac{\mathrm{d}A(w) }{|b-w|^2\log\frac{3}{|b-w|}}\geq\int_{D(b)}\frac{\mathrm{d}A(w) }{|b-w|^2\log\frac{3}{|b-w|}}.$$
Elementary calculations show that
$$\int_{D(b)}\frac{\mathrm{d}A(w) }{|b-w|^2\log\frac{3}{|b-w|}}=\frac{1}{\pi}\int_0^{2\pi}\mathrm{d}\theta\int_0^{1-|b|}\frac{1}{\rho\log \frac{3}{\rho}}\mathrm{d}\rho=\infty,$$
which shows that $\mathfrak{C}[g]\notin L^{\infty}(\ID)$.
\end{rem}
\subsection*{Proof of Theorem \ref{thm-Lp-BC*}}
Recall that
$$\mathfrak{J}_0[f](z)=\int_{\ID}\frac{{z}}{1-\bar{w}z}f(w)\mathrm{d}A(w) .$$
For $p>2$, assume that $f\in L^p(\ID)$. Applying H\"older's inequality for integrals, we have
\be\label{June-229-1}|\mathfrak{J}_0[f](z)|\leq\left(\int_{\ID}|f(w)|^p\mathrm{d}A(w) \right)^{\frac{1}{p}}\left(\int_{\ID}\frac{|z|^q}{|1-\bar{w}z|^q}\mathrm{d}A(w) \right)^{\frac{1}{q}},\ee
where $\frac{1}{p}+\frac{1}{q}=1$.

For $w=re^{it}\in\ID$, again by (\ref{2019-May-26-P}),  we have
\begin{align}\label{June-229-2}
\int_{\ID}\frac{|z|^q}{|1-\bar{w}z|^q}\mathrm{d}A(w) &=\frac{|z|^q}{\pi}\int_0^1r\mathrm{d}r\int_0^{2\pi}\frac{\mathrm{d}t}{|1-{z}re^{-it}|^q}\\ \nonumber
&=\sum\limits_{n=0}^{\infty}\left(\frac{\Gamma(n+\frac{q}{2})}{n!\Gamma(\frac{q}{2})}\right)^2\frac{|z|^{2n+q}}{n+1}\\ \nonumber
&= {}_2F_1\left[\frac{q}{2}, \frac{q}{2};2;|z|^{2}\right]|z|^q,
\end{align}
where $q=\frac{p}{p-1}\in[1, 2)$ and ${}_2F_1$ is the hypergeometric function.
Let
\be\label{June-31-1}M_q(|z|)={}_2F_1\left[\frac{q}{2}, \frac{q}{2};2;|z|^{2}\right]|z|^q.\ee
Then, $M_q(|z|)$ is an increasing function of $|z|$ and attains its maximum at $z=1$.
It follows from Lemma \Ref{DZ-2019-lemma-B} that
\be\label{June-229-3} M_q(1)=\frac{\Gamma(2-q)}{\Gamma^2\left(2-\frac{q}{2}\right)}.\ee
By using (\ref{June-229-1}), (\ref{June-229-2}) and (\ref{June-229-3}),
we have
$\|\mathfrak{J}_0[f]\|_{\infty}\leq \|f\|_{p} \left(\frac{\Gamma(2-q)}{\Gamma^2(2-\frac{q}{2})}\right)^{\frac{1}{q}},$
which implies that
\be\label{ZD-2020-Jan-10-1}\|\mathfrak{J}_0\|_{\infty}\leq\left(\frac{\Gamma(2-q)}{\Gamma^2(2-\frac{q}{2})}\right)^{\frac{1}{q}}.\ee

To show the equality of (\ref{ZD-2020-Jan-10-1}), fix $b\in\ID$, and consider the following function
$$g(w)=M_q^{-\frac{1}{p}}(|b|)\frac{\bar{b}}{1-w\bar{b}}\left|\frac{1-w\bar{b}}{\bar{b}}\right|^{\frac{p-2}{p-1}},$$
where $p>2$ and $q=\frac{p}{p-1}\in[1, 2)$.
Then
$$\|g\|_{p}^p=M_q^{-1}(|b|)\int_{\ID}\frac{|b|^q}{|1-\bar{w}b|^q}\mathrm{d}A(w) .$$
It follows from  (\ref{June-229-2}) and (\ref{June-31-1}) that
$\|g\|_{p}^p=1$, and thus $g\in L^p(\ID)$.

On the other hand, elementary calculations show that for $q=\frac{p}{p-1}<2$, we have
\begin{align*}
|\mathfrak{J}_0[g](b)|&=M_q^{-\frac{1}{p}}(|b|)\int_{\ID}\frac{|b|^q}{|1-{b}re^{-it}|^q}\mathrm{d}A(w) \\ \nonumber
&=M_q^{-\frac{1}{p}}(|b|)\sum\limits_{n=0}^{\infty}\left(\frac{\Gamma(n+\frac{q}{2})}{n!\Gamma(\frac{q}{2})}\right)^2\frac{|b|^{2n+q}}{n+1}\\ \nonumber
&=M_q^{1-\frac{1}{p}}(|b|).
\end{align*}
This shows that
$$\|\mathfrak{J}_0\|_{\infty}\geq\|\mathfrak{J}_0[g]\|_{\infty}\geq M_q^{1-\frac{1}{p}}(|b|).$$
Observe that, since $b$ can be arbitrarily close to $1$, by (\ref{June-229-3}), we have
\be\label{ZD-2020-Jan-10-2}\|\mathfrak{J}_0\|_{\infty}\geq \left(\frac{\Gamma(2-q)}{\Gamma^2(2-\frac{q}{2})}\right)^{\frac{1}{q}}.\ee
According to (\ref{ZD-2020-Jan-10-1}) and (\ref{ZD-2020-Jan-10-2}), we see that
$$\|\mathfrak{J}_0\|_{\infty}=\left(\frac{\Gamma(2-q)}{\Gamma^2\left(2-\frac{q}{2}\right)}\right)^{\frac{1}{q}},$$
and thus, (\ref{2019-June-30-1}) holds true.

If in particular $p=\infty$ (that is $q=1$), assume that $f\in L^\infty(\ID)$. Then using Lemma \ref{DZ-2019-lem-I1(w)}, we have
$$\|\mathfrak{J}_0[f]\|_\infty\leq\|f\|_\infty\int_{\ID}\frac{|z|}{|1-\bar{w}z|}\mathrm{d}A(w) \leq\frac{4}{\pi}\|f\|_\infty,$$
which shows that
\be\label{BZ-Dec-19}\|\mathfrak{J}_0\|_\infty\leq \frac{4}{\pi}.\ee

Fix $s\in\ID$, such that $s\neq0$. Let
$$f_s(w)=\frac{\bar{s}}{1-w\bar{s}}\left|\frac{1-w\bar{s}}{\bar{s}}\right|.$$
Then we get $\|f_s(w)\|_\infty=1$ and, by the proof of Lemma \ref{DZ-2019-lem-I1(w)}, we have
$$\mathfrak{J}_0[f_s](s)=\int_{\ID}\frac{|s|}{|1-\bar{w}s|}\mathrm{d}A(w) ={}_2F_1\bigg[\frac{1}{2}, \frac{1}{2}; 2; |s|^2\bigg]|s|.$$
This implies that
$$\|\mathfrak{J}_0[f_s]\|_\infty\geq{}_2F_1\bigg[\frac{1}{2}, \frac{1}{2}; 2; 1\bigg]=\frac{4}{\pi}.$$
Then
\be\label{BZ-Dec-19-1}\frac{4}{\pi}\leq\|\mathfrak{J}_0[f_s]\|_\infty\leq\|\mathfrak{J}_0\|_\infty.\ee
According to (\ref{BZ-Dec-19}) and (\ref{BZ-Dec-19-1}) we see that $\|\mathfrak{J}_0\|_\infty=\frac{4}{\pi}$.

This completes the proof of Theorem \ref{thm-Lp-BC*}. \qed
\begin{rem}\label{Remark-3.2-July}
For $p=2$, $\mathfrak{J}_0$ will not send $L^2(\ID)$ to $L^{\infty}(\ID)$. This can be seen as follows:
Let
$$g(w)=\frac{1}{(1-w)\log\frac{3}{|1-w|}},$$
where $w=\rho e^{it}\in\ID$.
Assume that $\xi=1-w=Re^{i\theta}$. Then
$$\int_{\ID}|g(w)|^2\mathrm{d}A(w) =\int_{\ID'}\frac{\mathrm{d}A(\xi)}{|\xi|^2\log^2\frac{3}{|\xi|}},$$
where $\ID'=\{\xi:|\xi-1|<1\}\subset \ID(0, 2):=\{\xi: |\xi|<2\}$. Therefore,
$$\int_{\ID}|g(w)|^2\mathrm{d}A(w) \leq\frac{1}{\pi}\int_0^{2\pi}\mathrm{d}\theta\int_0^2\frac{1}{R\log^2\frac{3}{R}}\mathrm{d}R=\frac{2}{\log\frac{3}{2}},$$
which shows that $g\in L^2(\ID)$.

Next, we are going to prove $\mathfrak{J}_0[g]\notin L^{\infty}(\ID)$. Once this is done, we then have $\mathfrak{J}_0$ doesn't send $L^p(\ID)$ to $L^{\infty}(\ID)$ for any $1\leq p\leq2$, because $L^2(\ID)\subseteq L^p(\ID)$.

For $z=r\in (0, 1)$ and $w=\rho e^{it}\in\ID$,
let
\begin{align*}
  G(t, \rho, r) & = \text{Re}\left(\frac{z}{1-\bar{w}z}g(w)\right)\\
  & =\frac{r (1+r \rho^2-\rho(1+r)\cos t)}{(1+r^2 \rho^2-2 r \rho \cos t) (1+\rho^2-2 \rho \cos t)  \log\frac{3}{\sqrt{1+\rho^2-2 \rho \cos t}}}.
\end{align*}
Then it is easy to see that $ G(t, \rho, r)>0$, for all $\rho, r\in(0, 1)$ and $t\in[0, 2\pi]$, since
$$1+r \rho^2-\rho(1+r)\cos t\geq (1-\rho)(1-r\rho)>0.$$
Therefore, we have
\begin{align}\label{april-29-1}
\nonumber|\mathfrak{J}_0[g](r)| & =\left|\int_{\ID}\frac{r}{1-\bar{w}r}g(w)\mathrm{d}A(w) \right| \\\nonumber
  & \geq \text{Re}\left(\int_{\ID}\frac{r}{1-\bar{w}r}g(w)\mathrm{d}A(w) \right)\\
  &=\frac{1}{\pi}\int_0^1\rho \mathrm{d}\rho\int_0^{2\pi}G(t, \rho, r)\mathrm{d}t.
\end{align}
Now, if the last integral of (\ref{april-29-1}) is infinity as $r\to 1$, then our problem is solved. To show this, by using Fatou's lemma (cf. \cite[Page 23]{Rudin}) and (\ref{april-29-1}), we have
$$  \varliminf_{r\to 1}  |\mathfrak{J}_0[g](r)|\ge \frac{1}{\pi}\int_{0}^1\rho \mathrm{d}\rho\int_0^{2\pi} \varliminf_{r\to 1} G(t,\rho,r) \mathrm{d}t, $$
where
$$\varliminf_{r\to 1}G(t,\rho,r)=\frac{1}{(1+\rho^2-2 \rho \cos t)\log\frac{3}{\sqrt{1+\rho^2-2 \rho \cos t}}}=\frac{1}{|1-w|^2\log\frac{3}{|1-w|}}.$$
Moreover,
\begin{align*}
  \frac{1}{\pi}\int_{0}^1\rho \mathrm{d}\rho\int_0^{2\pi} \varliminf_{r\to 1} G(t,\rho,r) dt& =\int_{\ID}\frac{1}{|1-w|^2\log\frac{3}{|1-w|}}\mathrm{d}A(w)  \\
  &= \int_{\ID'}\frac{1}{|\xi|^2\log\frac{3}{|\xi|}}\mathrm{d}A(\xi)\\
  &= \frac{1}{\pi}\int_{-\frac{\pi}{2}}^{\frac{\pi}{2}} \mathrm{d}\theta\int_0^{2\cos t} \frac{1}{R\log\frac{3}{R}}\mathrm{d}R=\infty.
\end{align*}
Based on the above discussions, we have $\mathfrak{J}_0[g]\notin L^{\infty}(\ID)$.
\end{rem}
\subsection*{Proof of Theorem \ref{2019-June-11-thm}}
For $p>2$, assume that $f\in L^p(\ID)$. It follows from H\"older's inequality for integrals that
\be\label{June-26-1}|\mathfrak{J}_0^*[f](z)|\leq\left(\int_{\ID}|f(w)|^p\mathrm{d}A(w) \right)^{\frac{1}{p}}\left(\int_{\ID}\frac{|w|^q}{|1-z\bar{w}|^q}\mathrm{d}A(w) \right)^{\frac{1}{q}},\ee
where $\frac{1}{p}+\frac{1}{q}=1$.

By using (\ref{2019-May-26-P}), for $w=re^{it}\in\ID$, we have
\begin{align}\label{June-18}
\int_{\ID}\frac{|w|^q}{|1-z\bar{w}|^q}\mathrm{d}A(w) &=\frac{1}{\pi}\int_0^1r^{q+1}\mathrm{d}r\int_0^{2\pi}\frac{\mathrm{d}t}{|1-zre^{-it}|^q}\\ \nonumber
&= 2\sum\limits_{n=0}^{\infty}\left(\frac{\Gamma(n+\frac{q}{2})}{n!\Gamma(\frac{q}{2})}\right)^2\frac{|z|^{2n}}{2n+q+2},
\end{align}
where $q=\frac{p}{p-1}\in[1, 2)$.

Let $N_q(|z|)=\int_{\ID}\frac{|w|^q}{|1-z\bar{w}|^q}\mathrm{d}A(w) $. Then according to (\ref{June-18}), we have $N_q(|z|)$ is an increasing function of $|z|$ and has its supremum at $z=1$.
Using the definition of the hypergeometric function, we have
\be\label{Mq1}N_q(1)=2\frac{{}_3F_2[1+\frac{q}{2}, \frac{q}{2},\frac{q}{2}; 1, 2+\frac{q}{2}; 1 ]}{2+q}.\ee
Moreover, it follows from Remark \ref{Rem-3.2} that $N_q(1)$ is finite.

Applying (\ref{June-26-1}) and (\ref{Mq1}), we see that
$$\|\mathfrak{J}_0^*[f]\|_{\infty}\leq \|f\|_{p}N_q^{1-\frac{1}{p}}(1).$$
Then
\be\label{June-29-3}
\|\mathfrak{J}_0^*\|_{\infty}\leq N_q^{1-\frac{1}{p}}(1).
\ee

To show the equality of (\ref{June-29-3}), fix $b\in\ID$, and consider the function
$$g(w)=\frac{w}{1-\bar{b}w}\left|\frac{1-\bar{b}w}{w}\right|^{\frac{p-2}{p-1}}N_q^{-\frac{1}{p}}(|b|),$$
where $p>2$.
Then for $q=\frac{p}{p-1}$, we have
$$\int_{\ID}|g(w)|^p\mathrm{d}A(w) =N_q^{-1}(|b|)\int_{\ID}\left|\frac{w}{1-\bar{b}w}\right|^q\mathrm{d}A(w) =1.$$
This shows that $\|g\|_p^p=1$.

On the other hand, elementary calculations show that
$$|\mathfrak{J}_0^*[g](b)|=N_q^{-\frac{1}{p}}(|b|)\int_{\ID}\left|\frac{w}{1-b\bar{w}}\right|^q\mathrm{d}A(w) =N_q^{1-\frac{1}{p}}(|b|).$$
Therefore, since $\sup\limits_{b\in\ID}N_q(|b|)=N_q(1)$, we have
$$\|\mathfrak{J}_0^*[g]\|_{\infty}=\sup\limits_{z\in\ID}|\mathfrak{J}_0^*[g](z)|\geq \lim\limits_{|b|\rightarrow1^-}|\mathfrak{J}_0^*[g](b)|= N_q^{1-\frac{1}{p}}(1).$$
This shows that
\be\label{ZD-2020-Jan-10-3}
\|\mathfrak{J}_0^*\|_{\infty}\geq N_q^{1-\frac{1}{p}}(1).
\ee
According to (\ref{Mq1}), (\ref{June-29-3}) and (\ref{ZD-2020-Jan-10-3}), we see that
(\ref{2019-June-11-1}) holds true.

If in particular $p=\infty$, assume that $f\in L^{\infty}(\ID)$. Then using Lemma \ref{DZ-2019-lem-I1(z)}, we see that
$$\|\mathfrak{J}_0^*[f]\|_{\infty}\leq\|f\|_{\infty}\sup\limits_{z\in\ID}\int_{\ID}\frac{|w|}{|1-\bar{w}z|}\mathrm{d}A(w) =\frac{1+2\alpha}{\pi}\|f\|_{\infty},$$
where $\alpha\approx 0.915966$ is the Catalan's constant.
Then
\be\label{ZD-2020-Jan-10-4}
\|\mathfrak{J}_0^*\|_{\infty}\leq\frac{1+2\alpha}{\pi}.
\ee

Fix $s\in\ID$, let
$$f_s(w)=\frac{1-\bar{w}s}{\bar{w}}\left|\frac{w}{1-\bar{w}s}\right|.$$
Then $\|f_s\|_{\infty}=1$ and by the proof of Lemma \ref{DZ-2019-lem-I1(z)}, we have
$$|\mathfrak{J}_0^*[f_s](s)| = 2\sum\limits_{n=0}^{\infty}\left(\frac{\Gamma(n+\frac{1}{2})}{n!\Gamma(\frac{1}{2})}\right)^2\frac{|s|^{2n}}{2n+3}.$$
This implies that $$\|\mathfrak{J}_0^*[f_s]\|_{\infty}\geq2\sum\limits_{n=0}^{\infty}\left(\frac{\Gamma(n+\frac{1}{2})}{n!\Gamma(\frac{1}{2})}\right)^2\frac{1}{2n+3}=\frac{1+2\alpha}{\pi}.$$
Then
\be\label{ZD-2020-Jan-10-5}
\|\mathfrak{J}_0^*\|_{\infty}\geq\|\mathfrak{J}_0^*[f_s]\|_{\infty}\geq\frac{1+2\alpha}{\pi}.
\ee

By (\ref{ZD-2020-Jan-10-4}) and (\ref{ZD-2020-Jan-10-5}), we have
$$\|\mathfrak{J}_0^*\|_{\infty}= \frac{1+2\alpha}{\pi}.$$

The proof of Theorem \ref{2019-June-11-thm} is complete. \qed

\begin{rem}\label{rem3.2}
For $p=2$, the image of $L^2(\ID)$ under $\mathfrak{J}_0^*$ is not contained in $L^{\infty}(\ID)$. In what follows, we will do more.
We construct a family of mappings $g_z$ continuous on $\overline{\ID}$, and that $\|g_z\|_2\leq\|g_1\|_2<\infty$. However,
$\mathfrak{J}_0^*[g_1]$ is not in $L^{\infty}(\ID)$.

For $z\in\ID$, set
$$g_z(w) = w(1-\bar{z} w)^{-1} \left(\log \frac{3}{|1-z\bar w|}\right)^{-1}.$$
Then for fixed $w\in\ID$, the mapping $h(z)= |g_z(w)|^2 $ is subharmonic, because direct computation leads to
$$\Delta h=16|w|^4A/B>0,$$
where $\Delta$ is the Laplace operator,
$C=\log|1-z\bar{w}|^2-\log9$, $A=(C+2)^2+2>0$ and $B=C^4|1-z\bar{w}|^2>0$.

Now, consider the following mapping
$$H(z)=\int_{\mathbb{D}} h(z) \mathrm{d}A(w) $$
which is subharmonic in $\ID$. Using the maximum principle for subharmonic functions, we get
\be\label{april-30-1} \|g_z(w)\|^2_{2}=\int_{\mathbb{D}} h(z) \mathrm{d}A(w) \le \max_{|z|=1} \int_{\mathbb{D}} h(z) \mathrm{d}A(w) =\int_{\mathbb{D}}h(1)\mathrm{d}A(w) . \ee

First, we prove $g_z\in L^2(\ID)$ as follows:
Let $\xi=1-w=Re^{i\theta}\in\ID'$. Then
$$\int_{\mathbb{D}}\frac{{|w|^2}}{|1-w|^2\log^2 \frac{3}{|1-w|}}\mathrm{d}A(w) \leq\int_{\ID(0, 2)}\frac{\mathrm{d}A(\xi)}{|\xi|^2\log^2\frac{3}{|\xi|}}=2\int_0^2\frac{1}{R\log^2\frac{3}{R}}\mathrm{d}R=\frac{2}{\log\frac{3}{2}}.$$
By using (\ref{april-30-1}), we see that $\|g_z(w)\|^2_{2}\leq\frac{2}{\log\frac{3}{2}}$, and thus,  $g_z\in L^2(\ID)$.

Second, we prove $\mathfrak{J}_0^*[g_1]\notin L^{\infty}(\ID)$ as follows: Following the proof of Remark \ref{Remark-3.2-July},
for $z=r\in(0, 1)$ and $w=\rho e^{it}\in\ID$, let
$$G_1(t, \rho, r) =\text{Re}\left(\frac{\bar{w}}{(1-r\bar{w})}g_1(w)\right).$$
Then
$$ G_1(t, \rho, r)=\frac{\rho^2}{r}G(t, \rho, r),$$
where $G(t, \rho, r)>0$ is the function given in Remark \ref{Remark-3.2-July}. Hence,
$G_1(t, \rho, r)>0$.

Moreover,
$$\varliminf_{r\to 1}G_1(t,\rho,r)=\rho^2G(t, \rho, 1)=\frac{|w|^2}{|1-w|^2\log\frac{3}{|1-w|}}.$$
Again, by Fatou's lemma, we see that for $\xi=1-w=Re^{i\theta}\in\ID'$,
\begin{align*}
  \varliminf_{r\to 1}|\mathfrak{J}_0^*[g_1](r)| &\geq\int_{\ID}\varliminf_{r\to 1}G_1(t,\rho,r)\mathrm{d}A(w) \\
  &=\int_{\mathbb{D}}\frac{{|w|^2}}{|1-w|^2\log \frac{3}{|1-w|}}\mathrm{d}A(w) \\
  &\geq\frac{1}{\pi}\int_{-\frac{\pi}{2}}^{\frac{\pi}{2}}\mathrm{d}\theta\int_0^{2\cos \theta}\frac{(1-R)^2}{R\log\frac{3}{R}}\mathrm{d}R.
\end{align*}
The divergence of the integral $\int_0^{2\cos \theta}\frac{(1-R)^2}{R\log(\frac{3}{R})}dR$ shows that $\varliminf_{r\to 1}|\mathfrak{J}_0^*[g_1](r)|$
is infinity, and thus, $\mathfrak{J}_0^*[g_1]\notin L^{\infty}(\ID)$. This shows that $\mathfrak{J}_0^*$ doesn't send $L^2(\ID)$ to
$L^{\infty}(\ID)$.
\end{rem}

\subsection*{Proof of Theorem \ref{thm-L2-BC}}
To prove Theorem \ref{thm-L2-BC}, following the proof of \cite[Theorem 5.2 and Corollary 5.3]{kalaj1}, it suffices to show that
$$\|\mathfrak{J}_0^*[P]\|_{2}^2\leq\frac{1}{2}\|P\|_{2}^2$$
whenever
$$P(z)=\sum\limits_{n=0}^{\infty}\sum\limits_{m=0}^{\infty}a_{m, n}z^m\bar{z}^n$$
is a polynomial of $z$ and $\bar{z}$, since such functions are dense in $L^2(\ID)$ and $\frac{1}{2}$ is the best constant.
In this case, only finitely many of the complex numbers $a_{m, n}$ are nonzero. It is evident that there exist radial functions
$f_d$, where $d$ is an integer number such that
$$P(z)=\sum\limits_{d=-\infty}^{\infty}g_d(z),$$
where $g_d(z)=f_d(r)e^{idt}$, $d=m-n$. Observe that $g_{d_1}$ and $g_{d_2}$ are orthogonal for $d_1\neq d_2$ in Hilbert space $L^2(\ID)$.

We will show that
\be\label{2020-May-2}\|\mathfrak{J}_0^*[P]\|_{2}^2\leq B\|P\|_{2}^2\ee
if and only if
\be\label{June-27-1}\sum\limits_{d=-\infty}^{\infty}\|\mathfrak{J}_0^*[g_d]\|_{2}^2\leq B\sum\limits_{d=-\infty}^{\infty}\|g_d\|_{2}^2,\ee
where $B=\frac{1}{2}$ is the best constant.

Suppose $w=re^{it}$ and $z=\rho e^{i\theta}\in\ID$.
By using the transform $\zeta=e^{it}$, we have
\begin{align*}
\mathfrak{J}_0^*[g_d](z)&=\int_{\ID}\frac{re^{-it}}{1-zre^{-it}}g_d(re^{it})\mathrm{d}A(w) \\ \nonumber
&=\frac{1}{\pi}\int_0^1r^2f_d(r)\mathrm{d}r\int_0^{2\pi}\frac{e^{i(d-1)t}}{1-zre^{-it}}\mathrm{d}t\\ \nonumber
&=\frac{1}{\pi}\int_0^1r^2f_d(r)\mathrm{d}r\int_{|\zeta|=1}\frac{\mathrm{d}\zeta}{i\zeta^{1-d}(\zeta-rz)}.
\end{align*}
Let $$\lambda_z(r)=\int_{|\zeta|=1}\frac{\mathrm{d}\zeta}{i\zeta^{1-d}(\zeta-rz)}.$$
Then by Cauchy residue theorem,
one has
$$\lambda_z(r)=\left\{
\begin{array}
{r@{\ }l}
2\pi(rz)^{d-1}, \ \  & \mbox{if}\ \  d\geq1;\\
\\
0 \ \ \ \ \ \ \ ,\ \ & \mbox{if}\ \  d<1.
\end{array}\right.$$
Now, we separate our discussions into two cases.
\begin{itemize}
  \item[{\bf Case 1.}]  Suppose $d\geq1$.
\end{itemize}
In this case,
$$\mathfrak{J}_0^*[g_d](z)=2A_dz^{d-1},$$
where
$$A_d=\int_0^1r^{d+1}f_d(r)\mathrm{d}r.$$
It is easy to see that $\mathfrak{J}_0^*[g_{d_1}]$ and $\mathfrak{J}_0^*[g_{d_2}]$ are orthogonal for any $d_1\neq d_2$, since
$$\int_{0}^{2\pi} e^{im\theta}\mathrm{d}\theta=0,\ \ \ \mbox{for any}\ \ \ m\neq0.$$
Therefore, we have $\mathfrak{J}_0^*[P](z)=\sum_{d=-\infty}^{\infty}\mathfrak{J}_0^*[g_d](z)$ and $\|\mathfrak{J}_0^*[P]\|_2^2=\sum_{d=-\infty}^{\infty}\|\mathfrak{J}_0^*[g_d]\|_{2}^2$.
Similarly, we can obtain $\|P\|_2^2=\sum_{d=-\infty}^{\infty}\|g_d\|_{2}^2.$
This shows that (\ref{2020-May-2}) and (\ref{June-27-1}) are equivalent.
Moreover, following the proof of the corresponding results in \cite[Theorem 5.2]{kalaj1} (see also \cite[Page 180]{Anderson}), to prove (\ref{June-27-1}), we only need to find the best constant $B$, such that
$\|\mathfrak{J}_0^*[g_d]\|_2^2\leq B\|g_d\|_2^2$, where $B=\frac{1}{2}$. In fact, we can only choose the function $P(w)=g_d(w)\in L^2(\ID)$, for fixed $d\in \mathbb{Z}$.

In what follows, we should find the best constant $B$ in (\ref{June-27-1}). Elementary calculations show that
\be\label{June-27-2}\int_{\ID}|\mathfrak{J}_0^*[g_{d}](z)|^2\mathrm{d}A(z) =\frac{4 |A_d|^2}{d}.\ee
On the other hand, we have
\be\label{June-27-3}\int_{\ID}|g_{d}(z)|^2\mathrm{d}A(z) =2\int_0^1\rho |f_d(\rho)|^2\mathrm{d}\rho.\ee
It follows from Proposition \ref{thm-Lp-BC} that there exists a constant $B$ such that
$$\int_{\ID}|\mathfrak{J}_0^*[g_d](z)|^2\mathrm{d}A(z) \leq B\int_{\ID}|g_{d}(z)|^2\mathrm{d}A(z) .$$
Applying (\ref{June-27-2}) and (\ref{June-27-3}), we have
$$\frac{4 |A_d|^2}{d}\leq2 B\int_0^1\rho |f_d(\rho)|^2\mathrm{d}\rho$$
that is
$$\left|\int_0^1r^{d+1}f_d(r)\mathrm{d}r\right|^2\leq\frac{dB}{2}\int_0^1r|f_d(r)|^2\mathrm{d}r.$$
To find the best constant $B$, using H\"older's inequality for integrals, we see that
$$\left|\int_0^1r^{d+1}f_d(r)\mathrm{d}r\right|^2\leq\frac{1}{2d+2}\int_0^1r|f_d(r)|^2\mathrm{d}r,$$
where the equality holds if $f_d(r)=Cr^{d}$ and $C$ is a constant.
This shows that
$$B=\frac{1}{d(d+1)},\ \ \ \mbox{for}\ \ \  d=1, 2, \cdots.$$ Thus,
$B=\frac{1}{2}$
is the best constant.

\begin{itemize}
  \item[{\bf Case 2.}]  Suppose $d<1$.
\end{itemize}
In this case,
$$\mathfrak{J}_0^*[g_d](z)=0.$$
Following the proof of Case 1, it is easy to see that in this case the best constant is $B=0$.

Based on the above discussions, we see that
$$B=\|\mathfrak{J}_0^*\|_{2}^2=\frac{1}{2}.$$

The proof of Theorem \ref{thm-L2-BC} is complete. \qed

\vspace*{5mm}
\noindent {\bf Acknowledgments}.
We would like to thank the anonymous referee for his/her helpful comments that have significant impact on this paper.

\vspace*{5mm}
\noindent {\bf Funding}.
The research of the authors were supported by NSFs of China (No. 11501220, 11971182), NSFs of Fujian Province (No. 2016J01020, 2019J0101)
and the Promotion Program for Young and Middle-aged Teachers in Science and Technology Research of Huaqiao University (ZQN-PY402).

\end{document}